\documentclass{amsart}
\usepackage{amsfonts}

\usepackage{graphicx}
\usepackage{amscd}
\usepackage{amsmath}

\newtheorem{theorem}{Theorem}
\theoremstyle{plain}

\newtheorem{corollary}{Corollary}

\newtheorem{lemma}{Lemma}

\newtheorem{proposition}{Proposition}

\numberwithin{equation}{section}

\begin{document}
\title{}
\author{}
\thanks{$^{\left( \ast \right) }$M\'{e}canique et Gravitation, Universit\'{e} de
Mons-Hainaut, Facult\'{e} des Sciences, 15 Avenue Maistriau, B-7000 Mons,
Belgium\\
Universit\`{a} degli Studi di Bergamo, Facolt\`{a} di Ingegneria, \noindent
Viale Marconi, 5, 24044 Dalmine \noindent (Bergamo) Italy \\
E-mail: Garattini@mi.infn.it}
\thanks{$^{\left( \ast \ast \right) }$ Memoria presentata da Marco Biroli, uno dei
XL.}
\maketitle

\begin{center}
{\large {\textbf{REMO GARATTINI}}}$^{\left( \ast \right) }$ \vspace{1cm}

\textrm{\textbf{{\large {GREEN'S FUNCTIONS AND ENERGY DECAY ON HOMOGENEOUS
SPACES}}}}$^{\left( \ast \ast \right) }$\vspace{1cm}
\end{center}

\noindent \textbf{Abstract.} We consider a homogeneous space $X=\left(
X,d,m\right) $ of dimension $\nu \geq 1$ and a local regular Dirichlet form
in $L^{2}\left( X,m\right) .$ We prove that if a Poincar\'{e} inequality
holds on every pseudo-ball $B\left( x,R\right) $ of $X$, with local
characteristic constant $c_{0}\left( x\right) $ and $c_{1}\left( r\right) $,
then a Green's function estimate from above and below is obtained. A
Saint-Venant-like principle is recovered in terms of the Energy's decay.%
\vspace{1cm}

\begin{center}
\textrm{\textbf{{\large {FUNZIONI DI GREEN E DECADIMENTO DELL'ENERGIA IN
SPAZI OMOGENEI}}}}
\end{center}

\vspace{.5cm}

\noindent \textbf{Sunto.} Si considera uno spazio omogeneo $X=\left(
X,d,m\right) $ di dimensione $\nu \geq 1$ e una forma di Dirichlet locale
regolare in $L^{2}\left( X,m\right) .$ Si dimostra che se una disuguaglianza
di Poincar\'{e} vale su ogni pseudo-sfera $B\left( x,R\right) $ di $X$, con
costanti caratteristiche locali $c_{0}\left( x\right) $ e $c_{1}\left(
r\right) $, allora si ricava una stima della funzione di Green da sopra e da
sotto. Il principio tipo Saint-Venant viene ottenuto in termini di
decadimento dell'Energia. \vspace{1cm}

\section{Introduction and Results}

We consider a \textit{connected, locally compact topological space }$X$. We
suppose that a \textit{distance} $d$ is defined on $X$ and we suppose that
the balls
\begin{equation*}
B\left( x,r\right) =\left\{ y\in X:d\left( x,y\right) <r\right\} ,\qquad r>0,
\end{equation*}
form a basis of open neighborhoods of $x\in X.$ Moreover, we suppose that a
(positive) Radon measure $m$ is given on $X$, with supp $m=X.$ The triple $%
\left( X,d,m\right) $ is assumed to satisfy the following property: there
exist some constants $0<R_{0}\leq +\infty ,\nu >0$ and a positive function $%
c_{0}\left( x\right) $ together with $c_{0}^{-1}\left( x\right) $ which
belong to $L_{loc}^{\infty }\left( X_{0}\right) $, where $X_{0}$ is a
relatively compact open subset of $X$, such that for every $x\in X$ and
every $0<r\leq R<R_{0}$%
\begin{equation}
0<c_{0}\left( x\right) \left( \frac{r}{R}\right) ^{\nu }m\left( B\left(
x,R\right) \right) \leq m\left( B\left( x,r\right) \right) .  \label{a1}
\end{equation}
Such a triple $\left( X,d,m\right) $ will be called a homogeneous space of
dimension $\nu .$ We point out, however, that a given exponent $\nu $
occurring in $\left( \ref{a1}\right) $ should be considered, more precisely,
as an upper bound of the ``homogeneous dimension'', hence we should better
call $\left( X,d,m\right) $ a homogeneous space of dimension less or equal
than $\nu .$ We further suppose that we are also given a \textit{strongly
local, regular, Dirichlet form a }in the Hilbert space $L^{2}\left(
X,m\right) $ - in the sense of M. Fukushima \cite{fuku}, - whose domain in $%
L^{2}\left( X,m\right) $ we shall denote by $\mathcal{D}\left[ a\right] $.
Furthermore, we shall restrict our study to Dirichlet forms of diffusion
type, that is to forms \textit{a} that have the following \textit{strong
local property:} $a\left( u,v\right) =0$ for every $u,v\in \mathcal{D}\left[
a\right] $ with $v$ constant on a neighborhood of supp $u.$ We recall that
the following integral representation of the form \textit{a }holds
\begin{equation*}
a\left( u,v\right) =\int_{X}\alpha \left( u,v\right) \left( dx\right)
,\qquad u,v\in \mathcal{D}\left[ a\right] ,
\end{equation*}
where $\alpha \left( u,v\right) $ is a uniquely defined signed Radon measure
on X, such that $\alpha \left( d,d\right) \leq m$, with $d\in \mathcal{D}%
_{loc}\left[ a\right] $: this last condition is fundamental for the
existence of cut off functions associated to the distance. Moreover, the
restriction of the measure $\alpha \left( u,v\right) $ to any open subset $%
\Omega $ of $X$, with $\overline{\Omega }\subset X_{0}$, depends only on the
restrictions of the functions $u,v$ to $\Omega .$ Therefore, the definition
of the measure $\alpha \left( u,v\right) $ can be unambiguously extended to
all $m$-measurable functions $u,v$ on $X$ that coincide $m-a.e.$ on every
compact subset of $\Omega $ with some functions of $\mathcal{D}\left[ a%
\right] .$ The space of all such functions will be denoted by $\mathcal{D}%
_{loc}\left[ a,\Omega \right] .$ Moreover we denote by $\mathcal{D}\left[
a,\Omega \right] $ the closure of $\mathcal{D}\left[ a\right] \cap C_{0}%
\left[ \Omega \right] $ in $\mathcal{D}\left[ a\right] $. The homogeneous
metric $d$ and the energy form $a$ associated to the \textit{energy measure }%
$\alpha $, both given on $X_{0}$, are then assumed to be mutually related by
the following basic assumption:

There exists a constant $k\geq 1$ such that $\forall $ $x\in X_{0}$, $%
\forall $ $R$ with $0<R<R_{0}$ the following Poincar\'{e} inequality holds
\cite{Biroli}:
\begin{equation}
\int\limits_{B\left( x,R\right) }\left| u-\bar{u}_{B\left( x,R\right)
}\right| ^{2}dm\leq c_{1}R^{2}\int\limits_{B\left( x,kR\right) }\alpha
\left( u,u\right) \left( dx\right)  \label{a2}
\end{equation}
for all $u\in D\left[ a,B\left( x,kR\right) \right] $, where
\begin{equation*}
\bar{u}_{B\left( x,R\right) }=\frac{1}{m\left( B\left( x,R\right) \right) }%
\int\limits_{B\left( x,R\right) }u\ dm,.
\end{equation*}
By assumption $\left( \ref{a2}\right) $, it can be shown the validity of the
following Sobolev type inequality of exponent $s$ for every $x\in X_{0}$ and
every $0<R<R_{0}$:
\begin{equation}
\left( \frac{1}{m\left( B\left( x,R\right) \right) }\int\limits_{B\left(
x,R\right) }\left| u\right| ^{s}dm\right) ^{\frac{1}{s}}\leq c_{1}R\left(
\int\limits_{B\left( x,R\right) }\alpha \left( u,u\right) \left( dx\right)
\right) ^{\frac{1}{2}},  \label{a3}
\end{equation}
where $u\in D\left[ a,B\left( x,kR\right) \right] $ and supp $u\subset
B\left( x,R\right) $. Let us consider the following simple generalization of
the Poincar\'{e} inequality
\begin{equation}
\int\limits_{B\left( x,R\right) }\left| u-\bar{u}_{B\left( x,R\right)
}\right| ^{2}dm\leq c_{1}^{2}\left( R\right) R^{2}\int\limits_{B\left(
x,kR\right) }\alpha \left( u,u\right) \left( dx\right)  \label{a5}
\end{equation}
where $c_{1}\left( r\right) $ is a decreasing function of $r$, then the
following Sobolev inequality of exponent $s$%
\begin{equation}
\left( \frac{1}{m\left( B\left( x,R\right) \right) }\int\limits_{B\left(
x,R\right) }\left| u\right| ^{s}dm\right) ^{\frac{1}{s}}\leq \tau ^{3}\left(
x\right) c_{1}\left( R\right) R\left( \int\limits_{B\left( x,R\right)
}\alpha \left( u,u\right) \left( dx\right) \right) ^{\frac{1}{2}}  \label{a6}
\end{equation}
can be proved, where we have defined $\tau \left( x\right) =\left( \underset{%
B\left( x,2R\right) }{\sup }c_{0}^{-1}\left( x\right) \right) ^{\frac{1}{2}}$%
. In this paper we will fix our attention on Green function estimates from
above and below by using an Harnack's inequality obtained in Ref.\cite{Remo}%
. Moreover, we will give the behaviour of the \textit{Energy decay} related
to the Green function under consideration. We begin here by recalling the
results given in \cite{Remo}.

\begin{theorem}[Harnack]
\label{tt1}Let $\left( \ref{a1}\right) $ , $\left( \ref{a5}\right) $, $%
\left( \ref{a6}\right) $ hold, and let $u$ be a non-negative solution of$\
a\left( u,v\right) =0.$ Let $\mathcal{O}$ be an open subset of $X_{0}$ and $%
u\in \mathcal{D}_{loc}\left[ \mathcal{O}\right] ,\ \forall v\in \mathcal{D}%
_{0}\left[ \mathcal{O}\right] $ with $B\left( x,r\right) \subset $ $\mathcal{%
O}$, then
\begin{equation*}
\underset{B_{\frac{1}{2}}}{ess\ \sup }\tilde{u}\leq \exp \gamma \mu \
\underset{B_{\frac{1}{2}}}{ess\ \inf }\tilde{u},
\end{equation*}
where $\tilde{u}$ is the function in $L^{2}\left( m,B\right) $ associated
with $u$, $\gamma $ $\equiv \gamma \left( \nu ,k\right) $, with $k$ a
positive constant and $\mu \left( x,r\right) =\tau ^{4}\left( x\right)
c_{1}\left( \frac{r}{2}\right) $. A standard consequence of the previous
Theorem is the following
\end{theorem}

\begin{corollary}
Suppose that
\begin{equation}
\int\limits_{r}^{R}e^{-\gamma \mu \left( x,\rho \right) }\frac{d\rho }{\rho }%
\rightarrow \infty \text{\qquad for\qquad }r\rightarrow 0  \label{aaa1}
\end{equation}
then the solution is continuous in the point under consideration. In
particular, if $\mu \left( x,\rho \right) \approx o\left( \log \log \frac{1}{%
\rho }\right) $, then there exists $c>0$ such that
\begin{equation}
\underset{B\left( x,r\right) }{osc}u\leq c\frac{\left( \log \frac{1}{R}%
\right) }{\left( \log \frac{1}{r}\right) }\underset{B\left( x,R\right) }{osc}%
u.
\end{equation}
\end{corollary}

Before presenting the main results we assume that $\partial B\left(
x,r\right) $ be connected and we prove the following

\begin{lemma}
\label{l0}Let $\bar{X}\equiv B\left( x,\frac{9}{4}r\right) -B\left( x,\frac{3%
}{4}r\right) $ be a connected set and let $l$ be equal to$\underset{x\in
B\left( x_{0},2R\right) }{\sup }c_{0}^{-1}\left( x\right) 16^{-\nu }$. Then
there exists a finite number $l$ of overlapping balls of radius $r/8$
joining two arbitrary points $x_{1}$ and $x_{l}$ of $\bar{X}$ which is at a
distance greater than $r/2$ from the origin.
\end{lemma}

\proof%
%
$\bar{X}$ can be covered with a finite set $\tilde{X}$ of balls of radius $%
r/16$. We can assume that the ball $B_{1}\equiv B\left( x_{1},r/16\right) $
is in the considered covering and we now assume that every ball in the
covering $\tilde{X}$ intersects $\bar{X}$, i.e. $B_{i}\cap \bar{X}\neq
\emptyset $ for every $i=1\ldots l$. Since $\bar{X}$ is a connected set,
there exists a second ball $B_{2}\equiv B\left( x_{2},r/16\right) \subset
\tilde{X}$ with $B_{2}\cap \bar{X}\neq \emptyset $ s.t. the closure of $%
B_{1} $ and the closure of $B_{2}$ do intersect, namely $\bar{B}_{1}\cap
\bar{B}_{2}\neq \emptyset $. Consider the set $X_{12}=\overline{B_{1}\cup
B_{2}}$. By hypothesis on $\bar{X}$ and $\tilde{X}$, there exists a ball $%
B_{3}\equiv B\left( x_{3},r/16\right) \subset \tilde{X}$, $B_{3}\cap \bar{X}%
\neq \emptyset $ s.t. $X_{12}\cap \bar{B}_{3}\neq \emptyset $. Now we can
consider the new set $X_{123}$ defined by $X_{123}=\overline{B_{1}\cup
B_{2}\cup B_{3}}$. By repeating the previous steps, we can say that there
exists a ball $B_{4}\equiv B\left( x_{4},r/16\right) \subset \tilde{X}$, $%
B_{4}\cap \bar{X}\neq \emptyset $ s.t. $X_{123}\cap \bar{B}_{4}\neq
\emptyset $. By iterating this procedure we can construct two sets $%
X_{1\ldots n}$ and $Y_{n+1\ldots l}$ s.t.$\underset{1\ldots n}{X}=\overline{%
\underset{i=1\ldots n}{\cup }B\left( x_{i},r/16\right) }$, $\underset{%
n+1\ldots l}{Y}=\overline{\underset{i=n+1\ldots l}{\cup }B\left(
x_{i},r/16\right) }$ with
\begin{equation*}
\underset{1\ldots n}{X}\cap \underset{n+1\ldots l}{Y}\neq \emptyset .
\end{equation*}
Indeed if $\underset{1\ldots n}{X}$ and $\underset{n+1\ldots l}{Y}$ were
s.t. $\underset{1\ldots n}{X}\cap \underset{n+1\ldots l}{Y}=\emptyset $ this
would mean that the space $\bar{X}$ is not connected against the hypothesis.
The chain of balls is obtained by joining the centers of the balls $B\left(
x_{i},r/16\right) $ forming the set
\begin{equation*}
\underset{1\ldots l}{X}=\overline{\underset{i=1\ldots l}{\cup }B\left(
x_{i},r/16\right) }
\end{equation*}
previously built which starts from $x_{1}$ and stops to $x_{l}$. To obtain
the overlapping of the balls is sufficient to consider $B\left(
x_{i},r/8\right) $ and the new sets
\begin{equation*}
\underset{1\ldots n}{\hat{X}}=\underset{i=1\ldots n}{\cup }B\left(
x_{i},r/8\right) ;\qquad \underset{n+1\ldots l}{\hat{Y}}=\underset{%
i=n+1\ldots l}{\cup }B\left( x_{i},r/8\right) ;
\end{equation*}
with
\begin{equation*}
\underset{1\ldots n}{\hat{X}}\cap \underset{n+1\ldots l}{\hat{Y}}\neq
\emptyset
\end{equation*}
for every $n$ and the Lemma is proved. Now we can state our main results

\begin{theorem}[Size of the Green function]
\label{tt2}$\forall B\left( x_{0},R\right) \subset B\left( x_{0},20R\right)
\subset \subset X_{0}$ and $\forall r\in \left( 0,\frac{R}{16}\right] $, the
following estimate holds for all $x\in \partial B\left( x_{0},r\right) $%
\begin{equation*}
\int\limits_{r}^{R}\frac{c_{0}\left( x\right) \exp \left( -l\gamma \mu
\left( x,s\right) \right) }{m\left( B\left( x,s\right) \right) }s^{2}\frac{ds%
}{s}\leq G_{B\left( x,R\right) }^{x}\left( y\right) \leq \int\limits_{r}^{R}%
\frac{\tau ^{6}\left( x\right) c_{1}^{2}\left( r\right) \exp \left( l\gamma
\mu \left( x,s\right) \right) }{m\left( B\left( x,s\right) \right) }s^{2}%
\frac{ds}{s},.
\end{equation*}
where $l=\underset{x\in B\left( x_{0},2R\right) }{\sup }c_{0}^{-1}\left(
x\right) 16^{-\nu }$ is a finite number of balls of radius $r/8$ covering $%
\partial B\left( x_{0},r\right) .$
\end{theorem}

\begin{theorem}[Saint-Venant-like principle]
\label{tt3}Let $u$ be a local solution in $X_{0}$ and $B\left(
x_{0},4R_{0}\right) \subset \subset X_{0}$ with $R_{0}=k^{2}R.$ Let
\begin{equation}
\psi \left( r\right) =\int\limits_{B\left( x_{0},r\right) }G_{B\left( x_{0},%
\frac{2r}{q}\right) }^{x_{0}}\alpha \left( u,u\right) \left( dx\right) .
\end{equation}
Then
\begin{equation*}
\psi \left( r\right) \leq cc_{1}^{4}\left( qr\right) \psi \left(
R_{0}\right) \exp \left( -\int\limits_{r}^{R}\exp \left( -2\gamma \mu \left(
x_{0},\rho \right) \frac{d\rho }{\rho }\right) \right)
\end{equation*}
\begin{equation}
+\exp \left( -\int\limits_{r}^{R}\exp \left( -5l\gamma \mu \left( x_{0},\rho
\right) \right) \frac{d\rho }{\rho }\right) \psi \left( R_{0}\right) .
\end{equation}
In particular, if $\gamma \mu \left( x,r\right) \leq o\left( \log \log \frac{%
1}{r}\right) $ then
\begin{equation}
\psi \left( r\right) \leq \left( \frac{\log \frac{1}{R}}{\log \frac{1}{r}}%
\right) \psi \left( R_{0}\right) ,
\end{equation}
with $k\geq 1.$
\end{theorem}

From the point of view of partial differential equations these results can
be applied to two important classes of operators on $\mathbb{R}^{n}$:

\begin{description}
\item[ a)]  \textit{Doubly Weighted uniformly elliptic operators} in
divergence form with measurable coefficients, whose coefficient matrix $%
A=\left( a_{ij}\right) $ satisfies
\begin{equation*}
w\left( x\right) \left| \xi \right| ^{2}\leq \left\langle A\xi ,\xi
\right\rangle \leq v\left( x\right) \left| \xi \right| ^{2}.
\end{equation*}
Here $\left\langle \cdot ,\cdot \right\rangle $ denotes the usual dot
product; $w$ and $v$ are weight functions, respectively belonging to $A_{2}$
and $D_{\infty }$ such that the following \textit{Poincar\'{e}} inequality
\begin{equation*}
\left( \frac{1}{\left| v\left( B\right) \right| }\int\limits_{B}\left|
f\left( x\right) -f_{B}\right| ^{2}v\ dx\right) ^{\frac{1}{2}}\leq cr\left(
\frac{1}{\left| w\left( B\right) \right| }\int\limits_{B}\left| \nabla
f\right| ^{2}w\ dx\right) ^{\frac{1}{2}}
\end{equation*}
holds.
\end{description}

\begin{description}
\item[b)]  \textit{Doubly Weighted H\"{o}rmander type operators }\cite
{franchi}$,$ whose form is $L=X_{k}^{\ast }\left( \alpha ^{hk}\left(
x\right) X_{h}\right) $ where $X_{h},h=1,\ldots ,m$ are smooth vector fields
in $\mathbb{R}^{n}$ that satisfy the H\"{o}rmander condition and $\alpha
=\left( \alpha ^{hk}\right) $ is any symmetric $m\times m$ matrix of
measurable functions on $\mathbb{R}^{n},$such that
\begin{equation*}
w\left( x\right) \sum\limits_{i}\left\langle X_{i},\xi \right\rangle
^{2}\leq \sum\limits_{i,j}\alpha _{ij}\left( x\right) \xi _{i}\xi _{j}\leq
v\left( x\right) \sum\limits_{i}\left\langle X_{i},\xi \right\rangle ^{2},
\end{equation*}
where $X_{i}\xi \left( x\right) =\left\langle X_{i},\nabla \xi \right\rangle
,\ i=1,\ldots ,m$, $\nabla \xi $ is the usual gradient of $\xi $ and $%
\left\langle \cdot ,\cdot \right\rangle $ denotes the usual inner product on
$\mathbb{R}^{n}.$ Then the following \textit{Poincar\'{e}} inequality for
vector fields
\begin{equation*}
\left( \frac{1}{\left| v\left( B\right) \right| }\int\limits_{B}\left|
f\left( x\right) -f_{B}\right| ^{2}v\ dx\right) ^{\frac{1}{2}}\leq cr\left(
\frac{1}{\left| w\left( B\right) \right| }\int\limits_{B}\left(
\sum\limits_{j}\left| \left\langle X_{j},\nabla f\left( x\right)
\right\rangle \right| ^{2}\right) ^{\frac{1}{2}}w\ dx\right) ^{\frac{1}{2}},
\end{equation*}
holds, with $w\in A_{2}$ and $v$ $\in D_{\infty }$.\textbf{\ }
\end{description}

\section{Estimates of the Green's functions and capacities of balls}

We define the \textit{Green} function $G_{\mathcal{O}}^{x}$ for the problem
\begin{equation}
\left\{
\begin{array}{c}
a\left( u,v\right) =\int\limits_{\mathcal{O}}fv\ m\left( dx\right) \\
u\in\mathcal{D}_{0}\left[ \mathcal{O}\right] ,\quad\forall v\in \mathcal{D}%
_{0}\left[ \mathcal{O}\right]
\end{array}
\right. ,  \label{21}
\end{equation}
$\mathcal{O}$ is a given ball $B\left( x_{0},R_{0}\right) \subset X_{0}$ and
$x\in\mathcal{O}$. The regularized Green function $G_{\rho,\mathcal{O}}^{x}$
associated with $\left( \ref{21}\right) $ is
\begin{equation}
\left\{
\begin{array}{c}
a\left( G_{\rho,\mathcal{O}}^{x},v\right) =%
-\hspace{-0.6cm}\int\limits_{B\left
( z,\rho\right) }%
%
vm\left( dx\right) \\
G_{\rho,\mathcal{O}}^{x}\in\mathcal{D}_{0}\left[ \mathcal{O}\right]
,\quad\forall v\in\mathcal{D}_{0}\left[ \mathcal{O}\right]
\end{array}
\right. ,
\end{equation}
where we have defined $%
-\hspace{-0.6cm}\int\limits_{B\left
( z,\rho\right) }%
%
=\left( \frac{1}{m\left( B\left( x,\rho\right) \right) }\right)
\int\limits_{B\left( x,\rho\right) }$, with $\rho>0$ and $B\left(
x,\rho\right) \subset\mathcal{O}.$ We define the \textit{capacity} of the
ball $B\left( x,r\right) $ with respect to the ball $B\left( x,dr\right) ,$ $%
d>1,$ relative to the form $a$, by setting
\begin{equation*}
cap\left( B\left( x,r\right) ,B\left( x,dr\right) \right) =\min\left\{
a\left( v,v\right) :v\in\mathcal{D}_{0}\left[ B\left( x,dr\right) \right]
,v\geq1\ m-a.e.\text{ on }B\left( x,r\right) \right\} .
\end{equation*}
By Sobolev-Poincar\'{e}'s inequality (jj), the minimum is achieved and the
unique minimizer $u\equiv u_{B\left( x,r\right) }$ is called the \textit{%
equilibrium potential} of $B\left( x,r\right) $ with respect to $B\left(
x,dr\right) ,$ relative to the form $a.$

\begin{theorem}
\label{t1} Let $G_{B\left( x,dr\right) }^{x}$ be the Green function of
problem $\left( \ref{21}\right) ,$ $\mathcal{O=}B\left( x,dr\right) ,$ with
singularity at $x,$ $d\geq 2,$ $B\left( x,4r\right) \subset \subset X_{0}.$
Suppose $\partial B\left( x,r\right) $ be connected, then the following
estimates hold: $\forall y\in \partial B\left( x,r\right) $
\begin{equation}
\frac{e^{-l\gamma \mu \left( x,r\right) }}{cap\left( B\left( x,r\right)
,B\left( x,dr\right) \right) }\leq G_{B\left( x,dr\right) }^{x}\left(
y\right) \leq \frac{e^{l\gamma \mu \left( x,r\right) }}{cap\left( B\left(
x,r\right) ,B\left( x,dr\right) \right) }  \label{22}
\end{equation}
and
\begin{equation}
c\frac{\left( d-1\right) ^{2}c_{0}\left( x\right) }{m\left( B\left(
x,r\right) \right) }r^{2}\leq \left( cap\left( B\left( x,r\right) ,B\left(
x,dr\right) \right) \right) ^{-1}\leq \frac{d^{2}r^{2}c_{1}^{2}\left(
r\right) \tau ^{6}\left( x\right) }{m\left( B\left( x,r\right) \right) },
\label{23}
\end{equation}
where $c_{1}\left( r\right) $ is the decreasing function of assumption $%
\left( \ref{a5}\right) $ $l=\underset{x\in B\left( x_{0},2R\right) }{\sup }%
c_{0}^{-1}\left( x\right) 16^{-\nu }$.
\end{theorem}

\proof%
%
Let us consider the cut-off function $\varphi $ of $B\left( x,r\right) $ in $%
B\left( x,\left( 1+\frac{d-1}{2}\right) r\right) $ as a test function, then
\begin{equation*}
cap\left( B\left( x,r\right) ,B\left( x,dr\right) \right) \leq
\int\limits_{B\left( x,dr\right) }\alpha \left( \varphi ,\varphi \right)
\left( dy\right)
\end{equation*}
\begin{equation}
\leq \frac{40d^{\nu }}{c_{0}\left( x\right) \left( d-1\right) ^{2}}\frac{%
m\left( B\left( x,r\right) \right) }{r^{2}}.  \label{24}
\end{equation}
There exists a positive Radon measure $\nu \equiv \nu _{B\left( x,r\right) }$
called the \textit{equilibrium measure} of $B\left( x,r\right) $ in $B\left(
x,dr\right) $ relative to the form $a$, such that
\begin{equation}
a\left( u_{B\left( x,r\right) },v\right) =\int\limits_{B\left( x,dr\right) }%
\tilde{v}\left( y\right) \nu _{B\left( x,r\right) }\left( dy\right)
\end{equation}
$\forall v\in \mathcal{D}_{0}\left[ B\left( x,dr\right) \right] .$ $\tilde{v}
$ is the q.c. version of $v,$ supp $\nu _{B\left( x,r\right) }\subset
\partial B\left( x,r\right) $ and
\begin{equation}
cap\left( B\left( x,r\right) ,B\left( x,dr\right) \right) =a\left(
u_{B\left( x,r\right) },u_{B\left( x,r\right) }\right) =\nu _{B\left(
x,r\right) }\left( \partial B\left( x,r\right) \right) .
\end{equation}
Since $u_{B\left( x,r\right) }\equiv 1,$ m-a.e. on $B\left( x,r\right) ,$ $%
\rho <\frac{r}{2},$ we have
\begin{equation*}
G_{\rho ,B\left( x,dr\right) }^{x}\in C\left( B\left( x,dr\right) -B\left( x,%
\frac{r}{2}\right) \right) \cap \mathcal{D}_{0}\left[ B\left( x,dr\right) %
\right] .
\end{equation*}
Then
\begin{equation}
a\left( u_{B\left( x,r\right) },G_{\rho ,B\left( x,dr\right) }^{x}\right)
=\int \alpha \left( u_{B\left( x,r\right) },G_{\rho ,B\left( x,dr\right)
}^{x}\right) \left( dy\right) =%
-\hspace{-0.6cm}\int\limits_{B\left
( z,\rho\right) }%
%
=1.
\end{equation}
But
\begin{equation*}
a\left( u_{B\left( x,r\right) },G_{\rho ,B\left( x,dr\right) }^{x}\right)
=\int\limits_{B\left( x,dr\right) }G_{\rho ,B\left( x,dr\right) }^{x}\left(
y\right) \nu _{B\left( x,r\right) }\left( dy\right)
\end{equation*}
\begin{equation}
=\int\limits_{\partial B\left( x,r\right) }G_{\rho ,B\left( x,dr\right)
}^{x}\left( y\right) \nu _{B\left( x,r\right) }\left( dy\right) =1
\end{equation}
and
\begin{equation*}
1=\int\limits_{\partial B\left( x,r\right) }G_{\rho ,B\left( x,dr\right)
}^{x}\left( y\right) \nu _{B\left( x,r\right) }\left( dy\right) \geq
\underset{\partial B\left( x,r\right) }{\inf }G_{\rho ,B\left( x,dr\right)
}^{x}\nu _{B\left( x,r\right) }\left( \partial B\left( x,r\right) \right)
\end{equation*}
\begin{equation}
\Longrightarrow 1\geq \underset{\partial B\left( x,r\right) }{\inf }G_{\rho
,B\left( x,dr\right) }^{x}cap\left( B\left( x,r\right) ,B\left( x,dr\right)
\right)
\end{equation}
Therefore
\begin{equation}
\frac{1}{cap\left( B\left( x,r\right) ,B\left( x,dr\right) \right) }\geq
\underset{\partial B\left( x,r\right) }{\inf }G_{\rho ,B\left( x,dr\right)
}^{x}\nu _{B\left( x,r\right) }  \label{25}
\end{equation}
But
\begin{equation*}
1=\int\limits_{\partial B\left( x,r\right) }G_{\rho ,B\left( x,dr\right)
}^{x}\left( y\right) \nu _{B\left( x,r\right) }\left( dy\right) \leq
\underset{\partial B\left( x,r\right) }{\sup }G_{\rho ,B\left( x,dr\right)
}^{x}\nu _{B\left( x,r\right) }\left( \partial B\left( x,r\right) \right) ,
\end{equation*}
this implies that
\begin{equation*}
1\leq \underset{\partial B\left( x,r\right) }{\sup }G_{\rho ,B\left(
x,dr\right) }^{x}cap\left( B\left( x,r\right) ,B\left( x,dr\right) \right)
\end{equation*}
that is
\begin{equation}
\frac{1}{cap\left( B\left( x,r\right) ,B\left( x,dr\right) \right) }\leq
\underset{\partial B\left( x,r\right) }{\sup }G_{\rho ,B\left( x,dr\right)
}^{x}.  \label{26}
\end{equation}
Collecting together $\left( \ref{25}\right) $ and $\left( \ref{26}\right) ,$
we obtain
\begin{equation}
\underset{\partial B\left( x,r\right) }{\inf }G_{\rho ,B\left( x,dr\right)
}^{x}\leq \frac{1}{cap\left( B\left( x,r\right) ,B\left( x,dr\right) \right)
}\leq \underset{\partial B\left( x,r\right) }{\sup }G_{\rho ,B\left(
x,dr\right) }^{x}.
\end{equation}
By Sobolev-Poincar\'{e} inequality, for $d\geq 2$, we have
\begin{equation*}
cap\left( B\left( x,r\right) ,B\left( x,dr\right) \right) =a\left(
u_{B\left( x,r\right) },u_{B\left( x,r\right) }\right)
\end{equation*}
\begin{equation}
\geq \frac{1}{d^{2}r^{2}c_{1}^{2}\left( r\right) \tau ^{6}\left( x\right) }%
\int\limits_{B\left( x,r\right) }u_{B\left( x,r\right) }^{2}m\left(
dx\right) =\frac{1}{d^{2}r^{2}c_{1}^{2}\left( r\right) \tau ^{6}\left(
x\right) }m\left( B\left( x,r\right) \right) ,
\end{equation}
which together with $\left( \ref{24}\right) $ shows that
\begin{equation*}
\frac{\left( d-1\right) ^{2}c_{0}\left( x\right) }{d^{\nu }40}\frac{r^{2}}{%
m\left( B\left( x,r\right) \right) }\leq \frac{1}{cap\left( B\left(
x,r\right) ,B\left( x,dr\right) \right) }
\end{equation*}
\begin{equation}
\leq \frac{d^{2}r^{2}c_{1}^{2}\left( r\right) \tau ^{6}\left( x\right) }{%
m\left( B\left( x,r\right) \right) }.
\end{equation}
By Lemma\ref{l0}, $\partial B\left( x,r\right) $ can be covered by a finite
number $l$ of overlapping balls of radius $r/8$ with distance greater than $%
r/4$ from $x$. This set has been denoted as
\begin{equation*}
\underset{1\ldots l}{\hat{X}}=\underset{i=1\ldots l}{\cup }B\left(
x_{i},r/8\right) .
\end{equation*}
Thus on each ball of $\underset{1\ldots l}{\hat{X}}$, a Harnack's inequality
holds
\begin{equation}
\underset{B_{i}}{\sup }G_{\rho ,B\left( x_{i},dr\right) }^{x_{i}}\leq
e^{\gamma \mu \left( x,r\right) }\underset{B_{i}}{\inf }G_{\rho ,B\left(
x_{i},dr\right) }^{x_{i}}.
\end{equation}
Let $u\left( x\right) =G_{\rho ,B\left( x,dr\right) }^{x}$. We begin with $%
B\left( x_{1},r/8\right) $ and $B\left( x_{2},r/8\right) $ both included in $%
\hat{X}$, then we obtain
\begin{equation*}
u\left( x_{1}\right) \leq \underset{B_{1}}{\sup }u\left( x\right) \leq
e^{\gamma \mu \left( x,r\right) }\underset{B_{1}}{\inf }u\left( x\right)
\leq e^{\gamma \mu \left( x,r\right) }u\left( \tilde{x}_{1}\right)
\end{equation*}
\begin{equation*}
\leq e^{\gamma \mu \left( x,r\right) }\underset{B_{2}}{\sup }u\left(
x\right) \leq \left( e^{\gamma \mu \left( x,r\right) }\right) ^{2}\underset{%
B_{2}}{\inf }u\left( x\right) \leq \left( e^{\gamma \mu \left( x,r\right)
}\right) ^{2}u\left( x_{2}\right) ,
\end{equation*}
where $\tilde{x}_{1}\in B\left( x_{1},r/8\right) \cap B\left(
x_{2},r/8\right) $. Let us consider the ball $B\left( x_{3},r/8\right)
\subset \underset{1\ldots l}{\hat{X}}$ as in Lemma\ref{l0}, then if $B\left(
x_{3},r/8\right) \cap B\left( x_{1},r/8\right) \neq \emptyset $ then
\begin{equation*}
u\left( x_{1}\right) \leq \left( e^{\gamma \mu \left( x,r\right) }\right)
^{2}u\left( x_{1}\right) ,
\end{equation*}
otherwise if $B\left( x_{3},r/8\right) \cap B\left( x_{1},r/8\right)
=\emptyset $ one gets
\begin{equation*}
u\left( x_{1}\right) \leq \left( e^{\gamma \mu \left( x,r\right) }\right)
^{3}u\left( x_{1}\right) .
\end{equation*}
By iterating the process to the $l$ balls of the chain, we get
\begin{equation*}
u\left( x_{1}\right) \leq \ldots \leq \left( e^{\gamma \mu \left( x,r\right)
}\right) ^{l}\underset{B_{l}}{\inf }u\left( x\right) \leq \left( e^{\gamma
\mu \left( x,r\right) }\right) ^{l}u\left( x_{l}\right) .
\end{equation*}
Then, collecting together the inequality chain and taking into account that $%
u\left( x\right) =G_{\rho ,B\left( x,dr\right) }^{x}$, we obtain
\begin{equation}
\underset{\partial B\left( x,r\right) }{\sup }G_{\rho ,B\left( x,dr\right)
}^{x}\leq e^{l\gamma \mu \left( x,r\right) }\underset{\partial B\left(
x,r\right) }{\inf }G_{\rho ,B\left( x,dr\right) }^{x}
\end{equation}
Therefore, by previous results we have
\begin{equation*}
\frac{1}{cap\left( B\left( x,r\right) ,B\left( x,dr\right) \right) }
\end{equation*}
\begin{equation*}
\leq \underset{\partial B\left( x,r\right) }{\sup }G_{\rho ,B\left(
x,dr\right) }^{x}\leq e^{l\gamma \mu \left( x,r\right) }\underset{\partial
B\left( x,r\right) }{\inf }G_{\rho ,B\left( x,dr\right) }^{x},
\end{equation*}
then
\begin{equation}
\frac{e^{-l\gamma \mu \left( x,r\right) }}{cap\left( B\left( x,r\right)
,B\left( x,dr\right) \right) }\leq \underset{\partial B\left( x,r\right) }{%
\inf }G_{\rho ,B\left( x,dr\right) }^{x}  \label{27}
\end{equation}
On the other hand
\begin{equation*}
\underset{\partial B\left( x,r\right) }{\inf }G_{\rho ,B\left( x,dr\right)
}^{x}\leq \frac{1}{cap\left( B\left( x,r\right) ,B\left( x,dr\right) \right)
}\leq \underset{\partial B\left( x,r\right) }{\sup }G_{\rho ,B\left(
x,dr\right) }^{x}
\end{equation*}
\begin{equation*}
\leq e^{l\gamma \mu \left( x,r\right) }\underset{\partial B\left( x,r\right)
}{\inf }G_{\rho ,B\left( x,dr\right) }^{x}\leq \frac{e^{l\gamma \mu \left(
x,r\right) }}{cap\left( B\left( x,r\right) ,B\left( x,dr\right) \right) },
\end{equation*}
then
\begin{equation}
\underset{\partial B\left( x,r\right) }{\sup }G_{\rho ,B\left( x,dr\right)
}^{x}\leq \frac{e^{l\gamma \mu \left( x,r\right) }}{cap\left( B\left(
x,r\right) ,B\left( x,dr\right) \right) }.  \label{28}
\end{equation}
Putting together $\left( \ref{27}\right) $ and $\left( \ref{28}\right) ,$ we
get
\begin{equation}
\frac{e^{-l\gamma \mu \left( x,r\right) }}{cap\left( B\left( x,r\right)
,B\left( x,dr\right) \right) }\leq G_{B\left( x,dr\right) }^{x}\left(
y\right) \leq \frac{e^{l\gamma \mu \left( x,r\right) }}{cap\left( B\left(
x,r\right) ,B\left( x,dr\right) \right) }.
\end{equation}
\textbf{Proof} \textbf{of Theorem} \textbf{\ref{tt2}}. Let $n\in \mathbb{N}$
be such that $2^{n}r<R<2^{n+1}r$ $\forall j=0,1,\ldots ,n$ and let $%
G_{j}^{x} $ be the Green function in $B\left( x,2^{j}r\right) $ with
singularity at $x$. Then by estimating from above and by Theorem \ref{t1} we
have
\begin{equation*}
G_{j}^{x}\left( y\right) \leq \tau ^{6}\left( x\right) \frac{c_{1}^{2}\left(
2^{j-1}r\right) \exp \left( l\gamma \mu \left( x,2^{j-1}r\right) \right) }{%
m\left( B\left( x,2^{j-1}r\right) \right) }d^{2}\left( 2^{j-1}r\right) ^{2},
\end{equation*}
with $y\in \partial B\left( x,2^{j-1}r\right) .$ We introduce the function
\begin{equation*}
u_{j}:=G_{j}^{x}-G_{j-1}^{x}\text{ \qquad in }B\left( x,2^{j-1}r\right) ,
\end{equation*}
which is a solution of $a\left( u_{j},v\right) =0.$ Indeed,
\begin{equation*}
a\left( u_{j},v\right) =a\left( G_{j}^{x},v\right) -a\left(
G_{j-1}^{x},v\right) =%
-\hspace{-0.6cm}\int\limits_{B\left
( x,\rho\right) }%
%
v\text{ }m\left( dx\right) -%
-\hspace{-0.6cm}\int\limits_{B\left
( x,\rho\right) }%
%
v\text{ }m\left( dx\right) =0,
\end{equation*}
with $u_{j}\in \mathcal{D}_{loc}\left[ B\left( x,2^{j-1}r\right) \right]
,\qquad \forall v\in \mathcal{D}_{0}\left[ B\left( x,2^{j-1}r\right) \right]
.$%
\begin{equation*}
\tilde{u}_{j}\left( y\right) -G_{j}^{x}\left( y\right) =-G_{j-1}^{x}\left(
y\right) =0\text{ \qquad q.e., }y\in \partial B\left( x,2^{j-1}r\right)
\end{equation*}
\begin{equation*}
\Longrightarrow \tilde{u}_{j}\left( y\right) \leq \frac{c_{j}\left(
r,x\right) }{m\left( B\left( x,2^{j-1}r\right) \right) }\left(
2^{j-1}r\right) ^{2}\qquad \text{q.e. on }\partial B\left( x,2^{j-1}r\right)
,\forall j=1,\ldots ,n
\end{equation*}
and $c_{j}\left( r,x\right) =\tau ^{6}\left( x\right) c_{1}^{2}\left(
2^{j-1}r\right) \exp \left( l\gamma \mu \left( x,2^{j-1}r\right) \right) $.
By the maximum principle
\begin{equation*}
u_{j}\left( y\right) \leq \frac{c_{j}\left( r,x\right) }{m\left( B\left(
x,2^{j-1}r\right) \right) }\left( 2^{j-1}r\right) ^{2}\qquad \text{m-a.e. in
}B\left( x,2^{j-1}r\right) ,j=1,\ldots ,n
\end{equation*}
if
\begin{equation*}
u:=G_{B\left( x,R\right) }^{x}-G_{n}^{x}\qquad \text{ in }B\left(
x,2^{n}r\right) ,
\end{equation*}
we find
\begin{equation*}
u\left( y\right) \leq \frac{c_{n}\left( r,x\right) }{m\left( B\left(
x,2^{n}r\right) \right) }\left( 2^{n}r\right) ^{2}\qquad \text{m-a.e. in }%
B\left( x,2^{n}r\right) .
\end{equation*}
This yields to
\begin{equation*}
G_{B\left( x,R\right) }^{x}\left( y\right) \leq u\left( y\right)
+\sum\limits_{j=1}^{n}u_{j}\left( y\right)
\end{equation*}
\begin{equation*}
\leq \frac{c_{n}\left( r,x\right) }{m\left( B\left( x,2^{n}r\right) \right) }%
\left( 2^{n}r\right) ^{2}+\sum\limits_{j=1}^{n}\frac{c_{j}\left( r,x\right)
}{m\left( B\left( x,2^{j-1}r\right) \right) }\left( 2^{j-1}r\right) ^{2}\leq
\sum\limits_{j=0}^{n}\frac{c_{j}\left( r,x\right) }{m\left( B\left(
x,2^{j-1}r\right) \right) }\left( 2^{j-1}r\right) ^{2}.
\end{equation*}
\begin{equation*}
\Longrightarrow G_{B\left( x,R\right) }^{x}\leq \sum\limits_{j=0}^{n}\tau
^{6}\left( x\right) \frac{c_{1}^{2}\left( 2^{j}r\right) \exp \left( l\gamma
\mu \left( x,2^{j}r\right) \right) }{m\left( B\left( x,2^{j}r\right) \right)
}d^{2}\left( 2^{j}r\right) ^{2}
\end{equation*}
\begin{equation*}
\Longrightarrow G_{B\left( x,R\right) }^{x}\left( y\right) \leq
\int\limits_{r}^{R}\tau ^{6}\left( x\right) \frac{c_{1}^{2}\left( r\right)
\exp \left( l\gamma \mu \left( x,s\right) \right) }{m\left( B\left(
x,s\right) \right) }s^{2}\frac{ds}{s}.
\end{equation*}
On the other hand, if we proceed to estimate the Green's function from
below, we have to consider the following initial inequality
\begin{equation*}
G_{j}^{x}\left( y\right) \geq \frac{c_{0}\left( x\right) \exp \left(
-l\gamma \mu \left( x,2^{j-1}r\right) \right) }{m\left( B\left(
x,2^{j-1}r\right) \right) }d^{2}\left( 2^{j-1}r\right) ^{2}.
\end{equation*}
By repeating the same steps of the estimate from above we arrive at
\begin{equation}
G_{B\left( x,R\right) }^{x}\left( y\right) \geq \int\limits_{r}^{R}\frac{%
c_{0}\left( x\right) \exp \left( -l\gamma \mu \left( x,s\right) \right) }{%
m\left( B\left( x,s\right) \right) }s^{2}\frac{ds}{s}
\end{equation}
and the desired estimate from above and below of the Green function becomes
\begin{equation*}
\int\limits_{r}^{R}\frac{c_{0}\left( x\right) \exp \left( -l\gamma \mu
\left( x,s\right) \right) }{m\left( B\left( x,s\right) \right) }s^{2}\frac{ds%
}{s}\leq G_{B\left( x,R\right) }^{x}\left( y\right) \leq
\int\limits_{r}^{R}\tau ^{6}\left( x\right) \frac{c_{1}^{2}\left( r\right)
\exp \left( l\gamma \mu \left( x,s\right) \right) }{m\left( B\left(
x,s\right) \right) }s^{2}\frac{ds}{s},
\end{equation*}
where we have taken the smallest radius on $c_{1}\left( r\right) $ because
of its decreasing property.

\section{Energy's decay}

We first prove a \textit{weighted Caccioppoli's inequality.}

\begin{proposition}
\label{pl}Let $v$ be a local solution in $B\left( x_{0},4r\right) ,$%
\begin{equation*}
\left\{
\begin{array}{c}
a\left( v,w\right) =0 \\
v\in \mathcal{D}_{loc}\left[ B\left( x_{0},4r\right) \right] ,\qquad \forall
w\in \mathcal{D}_{0}\left[ B\left( x_{0},4r\right) \right]
\end{array}
\right.
\end{equation*}
then
\begin{equation*}
\int\limits_{B\left( x_{0},qr\right) }G^{x_{0}}\alpha \left( v,v\right)
\left( dx\right) +\underset{B\left( x_{0},qr\right) }{\sup }v^{2}
\end{equation*}
\begin{equation}
\leq c_{q}\frac{s^{3l}\left( x_{0},r\right) c_{1}^{4}\left( qr\right)
s_{0}^{18}\left( x_{0}\right) }{m\left( B\left( x_{0},r\right) \right) }%
\int\limits_{B\left( x_{0},r\right) -B\left( x_{0},qr\right) }v^{2}m\left(
dx\right) ,
\end{equation}
where $c_{q}$ is a constant depending only by $q$ and where we have defined
\begin{equation}
s\left( x_{0},r\right) =\underset{z\in B\left( x_{0},r\right) -B\left(
x_{0},qr\right) }{\sup }e^{\gamma \mu \left( z,r\right) }  \label{30}
\end{equation}
and
\begin{equation}
s_{0}\left( x_{0}\right) =\underset{z\in B\left( x_{0},qr\right) }{\sup }%
\tau \left( z\right) ,  \label{30a}
\end{equation}
\end{proposition}

\proof%
%
Let $z\in B\left( x_{0},r\right) $ and $B\left( z,sr\right) \subset B\left(
z,tr\right) \subset B\left( x_{0},r\right) ,$ $s<t<1$ and let $\varphi$ be
the cut-off function of $B\left( z,sr\right) $ w.r.t. $B\left( z,tr\right) .$
We choose as test function $\varphi^{2}vG_{\rho}^{z},$ where $G_{\rho}^{z} $
denotes regularized Green function relative to $z$ and to the ball $B\left(
z,2r\right) .$ Since $\varphi,v,G_{\rho}^{z}\in\mathcal{D}_{loc}\left[
B\left( x_{0},2r\right) \right] \cap L^{\infty}\left( B\left(
x_{0},2r\right) ,m\right) $ we have
\begin{equation*}
0=a\left( v,\varphi^{2}vG_{\rho}^{z}\right) =\int\limits_{B\left(
z,tr\right) }\alpha\left( v,\varphi^{2}vG_{\rho}^{z}\right) \left( dx\right)
\end{equation*}
\begin{equation*}
=\int\limits_{B\left( z,tr\right) }\varphi^{2}G_{\rho}^{z}\alpha\left(
v,v\right) \left( dx\right) +\int\limits_{B\left( z,tr\right)
}\varphi^{2}v\alpha\left( v,G_{\rho}^{z}\right) \left( dx\right)
+\int\limits_{B\left( z,tr\right) }vG_{\rho}^{z}\alpha\left( v,\varphi
^{2}\right) \left( dx\right)
\end{equation*}
\begin{equation*}
=\int\limits_{B\left( z,tr\right) }\varphi^{2}G_{\rho}^{z}\alpha\left(
v,v\right) \left( dx\right) +\int\limits_{B\left( z,tr\right)
}\varphi^{2}v\alpha\left( v,G_{\rho}^{z}\right) \left( dx\right)
+\int\limits_{B\left( z,tr\right) }vG_{\rho}^{z}2\varphi\alpha\left(
v,\varphi\right) \left( dx\right)
\end{equation*}
\begin{equation*}
=\int\limits_{B\left( z,tr\right) }\varphi^{2}G_{\rho}^{z}\alpha\left(
v,v\right) \left( dx\right) +\frac{1}{2}\int\limits_{B\left( z,tr\right)
}\alpha\left( \varphi^{2}v^{2},G_{\rho}^{z}\right) \left( dx\right)
-\int\limits_{B\left( z,tr\right) }v^{2}\varphi\alpha\left( \varphi
,G_{\rho}^{z}\right) \left( dx\right)
\end{equation*}
\begin{equation}
+2\int\limits_{B\left( z,tr\right) }vG_{\rho}^{z}\varphi\alpha\left(
v,\varphi\right) \left( dx\right) .
\end{equation}
This implies that
\begin{equation*}
\int\limits_{B\left( z,tr\right) }\varphi^{2}G_{\rho}^{z}\alpha\left(
v,v\right) \left( dx\right) +\frac{1}{2}\int\limits_{B\left( z,tr\right)
}\alpha\left( \varphi^{2}v^{2},G_{\rho}^{z}\right) \left( dx\right)
\end{equation*}
\begin{equation*}
=\frac{1}{2}\int\limits_{B\left( z,tr\right) }2v^{2}\varphi\alpha\left(
\varphi,G_{\rho}^{z}\right) \left( dx\right) -2\int\limits_{B\left(
z,tr\right) }vG_{\rho}^{z}\varphi\alpha\left( v,\varphi\right) \left(
dx\right)
\end{equation*}
\begin{equation*}
\leq\int\limits_{B\left( z,tr\right) }\left| v^{2}\varphi\right| \left|
\alpha\left( \varphi,G_{\rho}^{z}\right) \right| \left( dx\right)
+2\int\limits_{B\left( z,tr\right) }\left| vG_{\rho}^{z}\varphi\right|
\left| \alpha\left( v,\varphi\right) \right| \left( dx\right)
\end{equation*}
\begin{equation*}
\leq\frac{1}{4}\int\limits_{B\left( z,tr\right)
}\varphi^{2}G_{\rho}^{z}\alpha\left( v,v\right) \left( dx\right)
+4\int\limits_{B\left( z,tr\right) }v^{2}G_{\rho}^{z}\alpha\left(
\varphi,\varphi\right) \left( dx\right)
\end{equation*}
\begin{equation}
+\frac{1}{4\varepsilon}\int\limits_{B\left( z,tr\right) -B\left( z,sr\right)
}v^{2}G_{\rho}^{z}\alpha\left( \varphi,\varphi\right) \left( dx\right)
+\varepsilon\int\limits_{B\left( z,tr\right) -B\left( z,sr\right)
}v^{2}\varphi^{2}\left( G_{\rho}^{z}\right) ^{-1}\alpha\left(
G_{\rho}^{z},G_{\rho}^{z}\right) \left( dx\right) .  \label{31}
\end{equation}

We will now estimate the last term at the r.h.s. of $\left( \ref{31}\right)
. $ Let $\sigma$ be the cut-off function of the annulus $B\left( z,tr\right)
-B\left( z,sr\right) $ w.r.t. the balls B$\left( z,\frac{s^{2}}{t}r\right) $
and B$\left( z,\left( 2t-s\right) r\right) $, then by applying Lemma (7.2)
of Ref.\cite{Biroli} with $f=\sigma\varphi v$, we find
\begin{equation*}
\int\limits_{B\left( z,tr\right) }\sigma^{2}v^{2}\varphi^{2}\left( G_{\rho
}^{z}\right) ^{-1}\alpha\left( G_{\rho}^{z},G_{\rho}^{z}\right) \left(
dx\right)
\end{equation*}
\begin{equation}
\leq4\int\limits_{B\left( z,tr\right) }v^{2}\varphi^{2}\left(
G_{\rho}^{z}\right) ^{2}\alpha\left( \sigma\varphi v,\sigma\varphi v\right)
\left( dx\right) .  \label{32}
\end{equation}
By means of Schwarz rule applied to $\alpha\left( \sigma\varphi
v,\sigma\varphi v\right) $, we have
\begin{equation*}
\int\limits_{B\left( z,tr\right) }\left( G_{\rho}^{z}\right)
^{2}\alpha\left( \sigma\varphi v,\sigma\varphi v\right) \left( dx\right)
\leq3\int\limits_{B\left( z,tr\right) -B\left( z,s^{\ast}r\right) }\left(
G_{\rho}^{z}\right) ^{2}\varphi^{2}\alpha\left( v,v\right) \left( dx\right)
\end{equation*}
\begin{equation}
+\frac{120}{\frac{s^{2}}{t^{2}}\left( t-s\right) ^{2}r^{2}}\int
\limits_{B\left( z,tr\right) -B\left( z,s^{\ast}r\right) }v^{2}\left(
G_{\rho}^{z}\right) ^{2}\text{ }m\left( dx\right) .  \label{33}
\end{equation}
Substituting $\left( \ref{33}\right) $ in $\left( \ref{32}\right) ,$ it
follows
\begin{equation*}
\int\limits_{B\left( z,tr\right) -B\left( z,s^{\ast}r\right)
}v^{2}\varphi^{2}\alpha\left( G_{\rho}^{z},G_{\rho}^{z}\right) \left(
dx\right) \leq\frac{480}{\frac{s^{2}}{t^{2}}\left( t-s\right) ^{2}r^{2}}%
\int\limits_{B\left( z,tr\right) -B\left( z,s^{\ast}r\right) }v^{2}\left(
G_{\rho}^{z}\right) ^{2}\text{ }m\left( dx\right)
\end{equation*}
\begin{equation}
+12\int\limits_{B\left( z,tr\right) -B\left( z,s^{\ast}r\right)
}\varphi^{2}\left( G_{\rho}^{z}\right) ^{2}\alpha\left( v,v\right) \left(
dx\right) .  \label{34}
\end{equation}
Let us consider
\begin{equation*}
\frac{\sup G_{\rho}^{z}}{\inf G_{\rho}^{z}},
\end{equation*}
where the $\sup$ and the $\inf$ are taken on $B\left( z,tr\right) -B\left(
z,s^{\ast}r\right) $ with $\rho<s^{\ast}r,$we have that
\begin{equation}
0<\frac{\sup G_{\rho}^{z}}{\inf G_{\rho}^{z}}\leq\tilde{c}\tau^{10}\left(
z\right) \exp\left( 2l\gamma\mu\left( z,r\right) \right) \left(
t/s^{\ast}\right) ^{\nu-2}c_{1}^{2}\left( r\right) .
\end{equation}
By multiplying inequality $\left( \ref{34}\right) $ by $\underset{z\in
B\left( z,tr\right) -B\left( z,s^{\ast}r\right) }{\sup\left(
G_{\rho}^{z}\right) ^{-1}}$ and taking account of estimates $\left( \ref{22}%
\right) $, we obtain
\begin{equation*}
\int\limits_{B\left( z,tr\right) -B\left( z,s^{\ast}r\right)
}v^{2}\varphi^{2}\left( G_{\rho}^{z}\right) ^{-1}\alpha\left(
G_{\rho}^{z},G_{\rho}^{z}\right) \left( dx\right)
\end{equation*}
\begin{equation*}
\leq\frac{\tau^{10}\left( z\right) 480\tilde{c}e^{2l\gamma\mu\left(
z,r\right) }c_{1}^{2}\left( r\right) }{\left( s/t\right) ^{2\nu-2}\left(
t-s\right) ^{2}r^{2}}\int\limits_{B\left( z,tr\right) -B\left( z,s^{\ast
}r\right) }v^{2}G_{\rho}^{z}\text{ }m\left( dx\right)
\end{equation*}
\begin{equation}
+12\tilde{c}\frac{\tau^{10}\left( z\right) e^{2l\gamma\mu\left( z,r\right)
}c_{1}^{2}\left( r\right) }{\left( s/t\right) ^{2\nu-4}}\int
\limits_{B\left( z,tr\right) }\varphi^{2}G_{\rho}^{z}\alpha\left( v,v\right)
\left( dx\right) .  \label{35}
\end{equation}
Putting $\left( \ref{34}\right) $ in $\left( \ref{31}\right) $, taking
account of the properties of $\varphi$ and choosing $\varepsilon\tilde{c}%
=\exp\left( -2l\gamma\mu\left( z,r\right) \right) \left( s/t\right)
^{2\nu-2}\tau^{-10}\left( z\right) /\left( 24c_{1}^{2}\left( r\right)
\right) ,$ we obtain
\begin{equation*}
\int\limits_{B\left( z,tr\right) }\varphi^{2}G_{\rho}^{z}\alpha\left(
v,v\right) \left( dx\right) +\frac{1}{2}\int\limits_{B\left( z,tr\right)
}\alpha\left( v^{2}\varphi^{2},G_{\rho}^{z}\right) \left( dx\right)
\end{equation*}
\begin{equation*}
\leq\frac{1}{4}\int\limits_{B\left( z,tr\right)
}\varphi^{2}G_{\rho}^{z}\alpha\left( v,v\right) \left( dx\right)
+4\int\limits_{B\left( z,tr\right) }v^{2}G_{\rho}^{z}\alpha\left(
\varphi,\varphi\right) \left( dx\right)
\end{equation*}
\begin{equation*}
+6\tilde{c}\tau^{10}\left( z\right) e^{2l\gamma\mu\left( z,r\right)
}c_{1}^{2}\left( r\right) \left( s/t\right) ^{2\nu-2}\int\limits_{B\left(
z,tr\right) -B\left( z,s^{\ast}r\right) }v^{2}G_{\rho}^{z}\alpha\left(
\varphi,\varphi\right) \left( dx\right)
\end{equation*}
\begin{equation*}
+\frac{20}{\left( t-s\right) ^{2}r^{2}}\int\limits_{B\left( z,tr\right)
-B\left( z,s^{\ast}r\right) }v^{2}G_{\rho}^{z}\text{ }m\left( dx\right) +%
\frac{1}{2}\frac{s^{2}}{t^{2}}\int\limits_{B\left( z,tr\right) }\varphi
^{2}G_{\rho}^{z}\alpha\left( v,v\right) \left( dx\right)
\end{equation*}
\begin{equation}
\leq\frac{c\tau^{10}\left( z\right) e^{2l\gamma\mu\left( z,r\right)
}c_{1}^{2}\left( r\right) \left( s/t\right) ^{2\nu-2}}{\left( t-s\right)
^{2}r^{2}}\int\limits_{B\left( z,tr\right) -B\left( z,s^{\ast}r\right)
}v^{2}G_{\rho}^{z}\text{ }m\left( dx\right) .
\end{equation}
Since $\varphi\equiv1$ on $B\left( z,sr\right) $ and by the definition of $%
G_{\rho}^{z}$ and recalling that $%
-\hspace{-0.6cm}\int\limits_{B\left
( z,\rho\right) }%
%
=\left( \frac{1}{m\left( B\left( z,\rho\right) \right) }\right)
\int\limits_{B\left( z,\rho\right) }$,

we have
\begin{equation*}
\int\limits_{B\left( z,tr\right) }\varphi ^{2}\left( G_{\rho }^{z}\right)
^{2}\alpha \left( v,v\right) \left( dx\right) +\frac{1}{2}%
-\hspace{-0.6cm}\int\limits_{B\left
( z,\rho\right) }%
%
v^{2}\text{ }m\left( dx\right)
\end{equation*}
\begin{equation}
\leq \frac{c\tau ^{10}\left( z\right) e^{2l\gamma \mu \left( z,r\right)
}c_{1}^{2}\left( r\right) \left( s/t\right) ^{2\nu -2}}{\left( t-s\right)
^{2}r^{2}}\int\limits_{B\left( z,tr\right) -B\left( z,s^{\ast }r\right)
}v^{2}G_{\rho }^{z}\text{ }m\left( dx\right) .
\end{equation}
Take the limit $\rho \rightarrow 0^{+}$, then $G_{\rho }^{z}\rightarrow
G^{z} $ uniformly in $B\left( z,tr\right) -B\left( z,s^{\ast }r\right) ,$ $%
\forall s^{\ast }$ fixed. By the Lebesgue theorem in Ref.\cite{Calderon}, we
obtain for m-a.e.
\begin{equation*}
\int\limits_{B\left( z,sr\right) }G^{z}\alpha \left( v,v\right) \left(
dx\right) +\frac{1}{2}\tilde{v}\left( z\right) ^{2}
\end{equation*}
\begin{equation}
\leq \frac{c\tau ^{10}\left( z\right) e^{2l\gamma \mu \left( z,r\right)
}c_{1}^{2}\left( r\right) \left( s/t\right) ^{2\nu -2}}{\left( t-s\right)
^{2}r^{2}}\int\limits_{B\left( z,tr\right) -B\left( z,s^{\ast }r\right)
}v^{2}G^{z}\text{ }m\left( dx\right) .
\end{equation}
We take the sup for $z\in B\left( x_{0},qr\right) $ by choosing $q\in \left(
0,\frac{1}{3}\right) ,\ s=\left[ 2q\left( 1-q\right) \right] ^{\frac{1}{2}%
},\ t=1-q.$ Then $B\left( z,tr\right) -B\left( z,s^{\ast }r\right) \subset
B\left( x_{0},r\right) -B\left( x_{0},qr\right) $ and $\ \forall z\in
B\left( x_{0},qr\right) $ we get
\begin{equation*}
\int\limits_{B\left( x_{0},qr\right) }G^{x_{0}}\alpha \left( v,v\right)
\left( dx\right) +\frac{1}{2}\underset{B\left( x_{0},qr\right) }{\sup }v^{2}
\end{equation*}
\begin{equation*}
\leq c_{q}\frac{s_{0}^{10}\left( x_{0}\right) s^{2l}\left( x_{0},r\right)
c_{1}^{2}\left( r\right) }{r^{2}}\underset{z\in B\left( x_{0},qr\right) }{%
\sup }\int\limits_{B\left( z,r\right) -B\left( z,qr\right) }v^{2}G^{z}\text{
}m\left( dx\right)
\end{equation*}
\begin{equation*}
\leq c_{q}\frac{s^{3l}\left( x_{0},r\right) c_{1}^{4}\left( qr\right)
s_{0}^{16}\left( x_{0}\right) }{m\left( B\left( z,qr\right) \right) }%
\int\limits_{B\left( x_{0},r\right) -B\left( x_{0},qr\right) }v^{2}\text{ }%
m\left( dx\right)
\end{equation*}
\begin{equation}
\leq c_{q}\frac{s^{3l}\left( x_{0},r\right) c_{1}^{4}\left( qr\right)
s_{0}^{18}\left( x_{0}\right) }{m\left( B\left( x_{0},r\right) \right) }%
\int\limits_{B\left( x_{0},r\right) -B\left( x_{0},qr\right) }v^{2}\text{ }%
m\left( dx\right)
\end{equation}
with $c_{q}=cq^{-\nu }\left( s/t\right) ^{2\nu -2}/\left( t-s\right) ^{2}$.

\textbf{Proof of Theorem \ref{tt3}}. Let us consider the test function $%
w=\left( u-k\right) G_{\rho }^{z}\varphi ,$ where $G_{\rho }^{z}$ is the
regularized Green function relative to $z$ and to the ball $B\left(
z,tr\right) ,$ $\varphi $ is the capacitory potential of $B\left(
z,sr\right) $ w.r.t. $B\left( z,tr\right) .$ $z\in B\left( x_{0},qr\right) ,$
$s<t<1,q$ to be fixed. $k$ is a constant. We have
\begin{equation*}
0=\int\limits_{B\left( z,tr\right) }\alpha \left( u,\left( u-k\right)
\varphi G_{\rho }^{z}\right) \left( dx\right)
\end{equation*}
\begin{equation*}
\int\limits_{B\left( z,tr\right) }\varphi G_{\rho }^{z}\alpha \left(
u,u\right) \left( dx\right) +\int\limits_{B\left( z,tr\right) }\left(
u-k\right) \varphi \alpha \left( u,G_{\rho }^{z}\right) \left( dx\right)
+\int\limits_{B\left( z,tr\right) }\left( u-k\right) G_{\rho }^{z}\alpha
\left( u,\varphi \right) \left( dx\right)
\end{equation*}
\begin{equation*}
=\int\limits_{B\left( z,tr\right) }\varphi G_{\rho }^{z}\alpha \left(
u,u\right) \left( dx\right) +\frac{1}{2}\int\limits_{B\left( z,tr\right)
}\alpha \left( \left( u-k\right) ^{2}\varphi ,G_{\rho }^{z}\right) \left(
dx\right)
\end{equation*}
\begin{equation}
+\int\limits_{B\left( z,tr\right) }\left( u-k\right) G_{\rho }^{z}\alpha
\left( u,\varphi \right) \left( dx\right) -\frac{1}{2}\int\limits_{B\left(
z,tr\right) }\left( u-k\right) ^{2}\alpha \left( \varphi ,G_{\rho
}^{z}\right) \left( dx\right) .
\end{equation}
Then for $\rho <r,$ we get
\begin{equation*}
\int\limits_{B\left( z,tr\right) }\varphi G_{\rho }^{z}\alpha \left(
u,u\right) \left( dx\right) +\frac{1}{2}%
-\hspace{-0.6cm}\int\limits_{B\left
( z,\rho\right) }%
%
\left( u-k\right) ^{2}m\left( dx\right)
\end{equation*}
\begin{equation*}
=\frac{1}{2}\int\limits_{B\left( z,tr\right) }\left( u-k\right) ^{2}\alpha
\left( \varphi ,G_{\rho }^{z}\right) \left( dx\right) -\int\limits_{B\left(
z,tr\right) }\left( u-k\right) G_{\rho }^{z}\alpha \left( u,\varphi \right)
\left( dx\right)
\end{equation*}
\begin{equation*}
=\frac{1}{2}\int\limits_{B\left( z,tr\right) }\alpha \left( \varphi ,G_{\rho
}^{z}\left( u-k\right) ^{2}\right) \left( dx\right) -2\int\limits_{B\left(
z,tr\right) }\left( u-k\right) G_{\rho }^{z}\alpha \left( u,\varphi \right)
\left( dx\right)
\end{equation*}
\begin{equation*}
=\frac{1}{2}\int\limits_{B\left( z,tr\right) }G_{\rho }^{z}\left( \tilde{u}%
-k\right) ^{2}d\nu _{B\left( z,sr\right) }-2\int\limits_{B\left( z,tr\right)
}\left( u-k\right) G_{\rho }^{z}\alpha \left( u,\varphi \right) \left(
dx\right)
\end{equation*}
\begin{equation*}
\leq \frac{1}{2}\underset{B\left( z,tr\right) }{\sup }\left( u-k\right)
^{2}+2\int\limits_{B\left( z,tr\right) }\left( u-k\right) G_{\rho
}^{z}\alpha \left( u,\varphi \right) \left( dx\right)
\end{equation*}
\begin{equation*}
\leq \frac{1}{2}\underset{B\left( z,tr\right) }{\sup }\left( u-k\right) ^{2}+%
\frac{1}{\eta }\int\limits_{B\left( z,tr\right) }G_{\rho }^{z}\alpha \left(
u,u\right) \left( dx\right) +\eta \int\limits_{B\left( z,tr\right) }\left(
u-k\right) ^{2}G_{\rho }^{z}\alpha \left( \varphi ,\varphi \right) \left(
dx\right)
\end{equation*}
\begin{equation*}
\leq \frac{1}{2}\underset{B\left( z,tr\right) }{\sup }\left( u-k\right) ^{2}+%
\frac{1}{\eta }\int\limits_{B\left( z,tr\right) }G_{\rho }^{z}\alpha \left(
u,u\right) \left( dx\right)
\end{equation*}
\begin{equation}
+\eta \underset{B\left( z,tr\right) }{\sup }\left( u-k\right) ^{2}\underset{%
B\left( z,tr\right) -B\left( z,sr\right) }{\sup }G_{\rho }^{z}\ cap\left(
B\left( z,sr\right) ,B\left( z,tr\right) \right) .
\end{equation}
Then, by the max principle and by Theorem $\left( \ref{t1}\right) $, we have
for arbitrary $\eta >0$%
\begin{equation*}
\int\limits_{B\left( z,tr\right) }\varphi G_{\rho }^{z}\alpha \left(
u,u\right) \left( dx\right) +\frac{1}{2}%
-\hspace{-0.6cm}\int\limits_{B\left
( z,\rho\right) }%
%
\left( u-k\right) ^{2}m\left( dx\right)
\end{equation*}
\begin{equation*}
\leq \frac{1}{2}\underset{B\left( z,tr\right) }{\sup }\left( u-k\right) ^{2}+%
\frac{1}{\eta }\int\limits_{B\left( z,tr\right) -B\left( z,sr\right)
}G_{\rho }^{z}\alpha \left( u,u\right) \left( dx\right)
\end{equation*}
\begin{equation}
+\eta \underset{B\left( z,tr\right) }{\sup }\left( u-k\right) ^{2}\underset{%
B\left( z,tr\right) -B\left( z,sr\right) }{\sup G_{\rho }^{z}}\ cap\left(
B\left( z,sr\right) ,B\left( z,tr\right) \right)
\end{equation}
Moreover, the application of Harnack's inequality on Green's function gives
\begin{equation*}
\underset{B\left( z,tr\right) -B\left( z,sr\right) }{\sup G_{\rho }^{z}}%
cap\left( B\left( z,sr\right) ,B\left( z,tr\right) \right)
\end{equation*}
\begin{equation*}
\leq e^{l\gamma \mu \left( z,r\right) }\underset{B\left( z,tr\right)
-B\left( z,sr\right) }{\inf G_{\rho }^{z}}cap\left( B\left( z,sr\right)
,B\left( z,tr\right) \right)
\end{equation*}
\begin{equation}
\leq c\tau ^{8}\left( z\right) c_{1}^{2}\left( r\right) e^{2l\gamma \mu
\left( z,r\right) },
\end{equation}
where $l=16^{-\nu }\sup\limits_{x\in B\left( x_{0},2R\right)
}c_{0}^{-1}\left( x\right) $. This implies that
\begin{equation*}
\int\limits_{B\left( z,tr\right) }\varphi G_{\rho }^{z}\alpha \left(
u,u\right) \left( dx\right) +\frac{1}{2}%
-\hspace{-0.6cm}\int\limits_{B\left
( z,\rho\right) }%
%
\left( u-k\right) ^{2}m\left( dx\right)
\end{equation*}
\begin{equation*}
\leq \frac{1}{2}\sup\limits_{B\left( z,tr\right) }\left( u-k\right)
^{2}+\eta cc\tau ^{8}\left( z\right) c_{1}^{2}\left( r\right)
\sup\limits_{z\in B\left( z,tr\right) -B\left( z,sr\right) }e^{2l\gamma \mu
\left( z,r\right) }\sup\limits_{B\left( z,tr\right) }\left( u-k\right) ^{2}
\end{equation*}
\begin{equation}
+\frac{1}{\eta }\int\limits_{B\left( z,tr\right) -B\left( z,sr\right)
}G_{\rho }^{z}\alpha \left( u,u\right) \left( dx\right)
\end{equation}
Then, passing to the limit as $\rho \rightarrow 0$ we obtain, by Lebesgue
theorem in \cite{Calderon}, for m-a.e.z
\begin{equation*}
\int\limits_{B\left( z,tr\right) }\varphi G_{B\left( x,tr\right) }^{z}\alpha
\left( u,u\right) \left( dx\right) +\frac{1}{2}\left( \tilde{u}\left(
z\right) -k\right) ^{2}
\end{equation*}
\begin{equation*}
\leq \left( \frac{1}{2}+\eta cc\tau ^{8}\left( z\right) c_{1}^{2}\left(
r\right) \sup\limits_{z\in B\left( z,tr\right) -B\left( z,sr\right)
}e^{2l\gamma \mu \left( z,r\right) }\right) \underset{B\left( z,tr\right) }{%
\sup }\left( u-k\right) ^{2}
\end{equation*}
\begin{equation}
+\frac{1}{\eta }\int\limits_{B\left( z,tr\right) -B\left( z,sr\right)
}G^{z}\alpha \left( u,u\right) \left( dx\right) .
\end{equation}
We now choose $t=1-q,s=2q,$ with $q\in \left( 0,\frac{1}{6}\right] .$ From
\begin{equation*}
\underset{B\left( z,tr\right) }{\sup }\left( u-k\right) ^{2}\leq \left(
1+2\eta cc\tau ^{8}\left( z\right) c_{1}^{2}\left( r\right)
\sup\limits_{z\in B\left( z,tr\right) -B\left( z,sr\right) }e^{2l\gamma \mu
\left( z,r\right) }\right) \underset{B\left( z,tr\right) }{\sup }\left(
u-k\right) ^{2}
\end{equation*}
\begin{equation}
+\frac{1}{\eta }\int\limits_{B\left( z,tr\right) -B\left( z,sr\right)
}G_{\rho }^{z}\alpha \left( u,u\right) \left( dx\right) ,
\end{equation}
we take the supremum for $z\in B\left( x_{0},qr\right) $. Then from Theorem $%
\left( \ref{t1}\right) $ we have
\begin{equation*}
\underset{B\left( x_{0},qr\right) }{\sup }\left( u-k\right) ^{2}\leq \left(
1+2\eta cs_{0}^{8}\left( x_{0}\right) c_{1}^{2}\left( r\right) s^{2l}\left(
x_{0},r\right) \right) \underset{B\left( x_{0},r\right) }{\sup }\left(
u-k\right) ^{2}
\end{equation*}
\begin{equation*}
+\frac{1}{\eta }\left( \frac{s^{2l}\left( x_{0},r\right) }{cap\left( B\left(
x_{0},r\right) ,B\left( x_{0},qr\right) \right) }\right)
\int\limits_{B\left( x_{0},r\right) -B\left( x_{0},qr\right) }\alpha \left(
u,u\right) \left( dx\right)
\end{equation*}
\begin{equation*}
\leq \left( 1+2\eta cs_{0}^{8}\left( x_{0}\right) c_{1}^{2}\left( r\right)
s^{2l}\left( x_{0},r\right) \right) \underset{B\left( x_{0},r\right) }{\sup }%
\left( u-k\right) ^{2}
\end{equation*}
\begin{equation}
+\frac{s^{2l}\left( x_{0},r\right) }{\eta }\int\limits_{B\left(
x_{0},r\right) -B\left( x_{0},qr\right) }e^{l\gamma \mu \left(
x_{0},r\right) }G_{B\left( x_{0},2r\right) }^{z}\alpha \left( u,u\right)
\left( dx\right) ,  \label{tl1}
\end{equation}
with $s\left( x_{0},r\right) $ defined as in Eq.$\left( \ref{30}\right) $.
From Proposition $\left( \ref{pl}\right) $, it follows that
\begin{equation*}
\int\limits_{B\left( x_{0},qr\right) }G^{x_{0}}\alpha \left( v,v\right)
\left( dx\right) +\frac{1}{2}\underset{B\left( x_{0},qr\right) }{\sup }v^{2}
\end{equation*}
\begin{equation*}
\leq c_{q}\frac{s^{3l}\left( x_{0},r\right) c_{1}^{4}\left( qr\right)
s_{0}^{18}\left( x_{0}\right) }{m\left( B\left( x_{0},r\right) \right) }%
\int\limits_{B\left( x_{0},r\right) -B\left( x_{0},qr\right) }v^{2}\text{ }%
m\left( dx\right)
\end{equation*}
\begin{equation}
\leq c_{q}c_{1}^{4}\left( qr\right) s^{3l}\left( x_{0},r\right)
s_{0}^{18}\left( x_{0}\right) \sup\limits_{B\left( x_{0},r\right) }v^{2},
\end{equation}
where $s_{0}\left( x_{0}\right) $ has been defined in Eq.$\left( \ref{30a}%
\right) $. Thus
\begin{equation*}
\int\limits_{B\left( x_{0},qr\right) }G_{B\left( x_{0},2r\right)
}^{x_{0}}\alpha \left( u,u\right) \left( dx\right)
\end{equation*}
\begin{equation}
\leq c_{q}c_{1}^{4}\left( qr\right) s^{3l}\left( x_{0},r\right)
s_{0}^{18}\left( x_{0}\right) \sup\limits_{B\left( x_{0},r\right) }\left(
u-k\right) ^{2}.
\end{equation}
Then, by $\left( \ref{tl1}\right) $ we obtain
\begin{equation*}
\int\limits_{B\left( x_{0},qr\right) }G_{B\left( x_{0},2r\right)
}^{x_{0}}\alpha \left( u,u\right) \left( dx\right) +\sup\limits_{B\left(
x_{0},qr\right) }\left( u-k\right) ^{2}
\end{equation*}
\begin{equation*}
\leq \left[ c_{q}c_{1}^{4}\left( qr\right) s^{3l}\left( x_{0},r\right)
s_{0}^{18}\left( x_{0}\right) +1+2\eta cs_{0}^{8}\left( x_{0}\right)
c_{1}^{2}\left( r\right) s^{2l}\left( x_{0},r\right) \right] \underset{%
B\left( x_{0},r\right) }{\sup }\left( u-k\right) ^{2}
\end{equation*}
\begin{equation}
+\frac{s^{3l}\left( x_{0},r\right) }{\eta }\int\limits_{B\left(
x_{0},r\right) -B\left( x_{0},qr\right) }G_{B\left( x_{0},2r\right)
}^{x_{0}}\alpha \left( u,u\right) \left( dx\right) .
\end{equation}
By ``hole filling'' after having multiplied by $\eta $, we obtain
\begin{equation*}
\left( s^{3l}\left( x_{0},r\right) +\eta \right) \int\limits_{B\left(
x_{0},qr\right) }G_{B\left( x_{0},2r\right) }^{x_{0}}\alpha \left(
u,u\right) \left( dx\right) +\eta \underset{B\left( x_{0},qr\right) }{\sup }%
\left| u-k\right| ^{2}
\end{equation*}
\begin{equation*}
\leq \eta \underset{B\left( x_{0},r\right) }{\sup }\left( u-k\right) ^{2}%
\left[ c_{q}c_{1}^{4}\left( qr\right) s^{3l}\left( x_{0},r\right)
s_{0}^{18}\left( x_{0}\right) +1+2\eta cs_{0}^{8}\left( x_{0}\right)
c_{1}^{2}\left( r\right) s^{2l}\left( x_{0},r\right) \right]
\end{equation*}
\begin{equation}
+s^{3l}\left( x_{0},r\right) \int\limits_{B\left( x_{0},r\right) }G_{B\left(
x_{0},2r\right) }^{x_{0}}\alpha \left( u,u\right) \left( dx\right) .
\label{tl2}
\end{equation}
We now study the last term at the right hand side of $\left( \ref{tl2}%
\right) $%
\begin{equation*}
\int\limits_{B\left( x_{0},r\right) }G_{B\left( x_{0},2r\right)
}^{x_{0}}\alpha \left( u,u\right) \left( dx\right) =\int\limits_{B\left(
x_{0},r\right) }G_{B\left( x_{0},2q^{-1}r\right) }^{x_{0}}\alpha \left(
u,u\right) \left( dx\right)
\end{equation*}
\begin{equation}
-\int\limits_{B\left( x_{0},r\right) }\left( G_{B\left(
x_{0},2q^{-1}r\right) }^{x_{0}}-G_{B\left( x_{0},2r\right) }^{x_{0}}\right)
\alpha \left( u,u\right) \left( dx\right) .
\end{equation}
Here we have taken into account that
\begin{equation*}
F=G_{B\left( x_{0},2q^{-1}r\right) }^{x_{0}}-G_{B\left( x_{0},2r\right)
}^{x_{0}}
\end{equation*}
is a solution of the problem
\begin{equation*}
a\left( F,v\right) =0,\qquad \forall v\in \mathcal{D}\left[ B\left(
x_{0},2r\right) \right] ,
\end{equation*}
therefore by the maximum principle and Theorem $\left( \ref{t1}\right) $%
\begin{equation*}
\underset{B\left( x_{0},r\right) }{\inf }F\geq \underset{\partial B\left(
x_{0},2r\right) }{\inf }\tilde{F}=\underset{\partial B\left( x_{0},2r\right)
}{\inf }G_{B\left( x_{0},2q^{-1}r\right) }^{x_{0}}
\end{equation*}
\begin{equation}
\geq c\exp \left( -l\gamma \mu \left( x_{0},r\right) \right) \left( \frac{%
c_{0}\left( x_{0}\right) r^{2}}{m\left( B\left( x_{0},r\right) \right) }%
\right) .
\end{equation}
Therefore, by Poincar\'{e} inequality, we also have for arbitrary $\bar{q}%
\in \left( 0,1\right) $%
\begin{equation*}
\int\limits_{B\left( x_{0},r\right) }G_{B\left( x_{0},2r\right)
}^{x_{0}}\alpha \left( u,u\right) \left( dx\right) \leq \int\limits_{B\left(
x_{0},r\right) }G_{B\left( x_{0},2q^{-1}r\right) }^{x_{0}}\alpha \left(
u,u\right) \left( dx\right)
\end{equation*}
\begin{equation*}
-c\exp \left( -l\gamma \mu \left( x_{0},r\right) \right) \left( \frac{%
c_{0}\left( x_{0}\right) r^{2}}{m\left( B\left( x_{0},r\right) \right) }%
\right) \int\limits_{B\left( x_{0},r\right) }\alpha \left( u,u\right) \left(
dx\right)
\end{equation*}
\begin{equation*}
\leq \int\limits_{B\left( x_{0},r\right) }G_{B\left( x_{0},2q^{-1}r\right)
}^{x_{0}}\alpha \left( u,u\right) \left( dx\right)
\end{equation*}
\begin{equation}
-c\left( \frac{c_{0}\left( x_{0}\right) \exp \left( -\gamma \mu \left(
x_{0},r\right) \right) r^{2}}{m\left( B\left( x_{0},r\right) \right)
c_{1}\left( \kappa ^{-1}\bar{q}r\right) \left( \kappa ^{-1}\bar{q}r\right)
^{2}}\right) \int\limits_{B\left( x_{0},\kappa ^{-1}\bar{q}r\right) }\left|
u-\bar{u}\right| ^{2}m\left( dx\right) ,
\end{equation}
where $\bar{u}$ denotes the average of $u$ on $B\left( x_{0},\kappa ^{-1}%
\bar{q}r\right) .$ By choosing $\bar{q}$ such that $\kappa ^{-1}\bar{q}=q,$
we find
\begin{equation*}
\int\limits_{B\left( x_{0},r\right) }G_{B\left( x_{0},2r\right)
}^{x_{0}}\alpha \left( u,u\right) \left( dx\right) \leq \int\limits_{B\left(
x_{0},r\right) }G_{B\left( x_{0},2q^{-1}r\right) }^{x_{0}}\alpha \left(
u,u\right) \left( dx\right)
\end{equation*}
\begin{equation}
-c\left( \frac{c_{0}\left( x_{0}\right) \exp \left( -\gamma \mu \left(
x_{0},r\right) \right) }{m\left( B\left( x_{0},r\right) \right) }\right)
\frac{1}{c_{1}\left( qr\right) }\underset{B\left( x_{0},qr\right) }{\sup }%
\left| u-\bar{u}\right| ^{2}m\left( B\left( x_{0},qr\right) \right) ,
\end{equation}
while taking the doubling property of $m$ into account,
\begin{equation*}
\frac{m\left( B\left( x_{0},qr\right) \right) }{m\left( B\left(
x_{0},r\right) \right) }\geq c_{0}\left( x_{0}\right) q^{\nu },
\end{equation*}
we obtain that
\begin{equation*}
\int\limits_{B\left( x_{0},r\right) }G_{B\left( x_{0},2r\right)
}^{x_{0}}\alpha \left( u,u\right) \left( dx\right) \leq \int\limits_{B\left(
x_{0},r\right) }G_{B\left( x_{0},2q^{-1}r\right) }^{x_{0}}\alpha \left(
u,u\right) \left( dx\right)
\end{equation*}
\begin{equation}
-cq^{\nu -2}\left( \exp \left( -l\gamma \mu \left( x_{0},r\right) \right)
\right) \frac{c_{0}^{2}\left( x_{0}\right) }{c_{1}\left( qr\right) }%
\underset{B\left( x_{0},qr\right) }{\sup }\left| u-\bar{u}\right| ^{2}.
\label{tl3}
\end{equation}
Taking into account $\left( \ref{tl3}\right) $ and choosing $\bar{u}=k$ in $%
\left( \ref{tl2}\right) ,$ we obtain
\begin{equation*}
\left( s^{3l}\left( x_{0},r\right) +\eta \right) \int\limits_{B\left(
x_{0},qr\right) }G_{B\left( x_{0},2r\right) }^{x_{0}}\alpha \left(
u,u\right) \left( dx\right)
\end{equation*}
\begin{equation}
+\left( \eta +cq^{\nu -2}\left( \exp \left( -l\gamma \mu \left(
x_{0},r\right) \right) \right) \frac{c_{0}^{2}\left( x_{0}\right) }{%
c_{1}\left( qr\right) }\right) \underset{B\left( x_{0},qr\right) }{\sup }%
\left| u-\bar{u}\right| ^{2}
\end{equation}
\begin{equation*}
\leq \eta \underset{B\left( x_{0},r\right) }{\sup }\left| u-\bar{u}\right|
^{2}\left[ c_{q}c_{1}^{4}\left( qr\right) s^{3l}\left( x_{0},r\right)
s_{0}^{18}\left( x_{0}\right) +1+2\eta cs_{0}^{8}\left( x_{0}\right)
c_{1}^{2}\left( r\right) s^{2l}\left( x_{0},r\right) \right]
\end{equation*}
\begin{equation}
+s^{3l}\left( x_{0},r\right) \int\limits_{B\left( x_{0},r\right) }G_{B\left(
x_{0},2r\right) }^{x_{0}}\alpha \left( u,u\right) \left( dx\right) .
\end{equation}
Since
\begin{equation}
s\left( x_{0},r\right) =\underset{z\in B\left( x_{0},r\right) -B\left(
x_{0},qr\right) }{\sup e^{\gamma \mu \left( x_{0},r\right) }}=\exp \bar{s}%
\left( x_{0},r\right)
\end{equation}
with $\bar{s}\left( x_{0},r\right) =\underset{z\in B\left( x_{0},r\right)
-B\left( x_{0},qr\right) }{\sup \gamma \mu \left( x_{0},r\right) }$, we now
choose $\eta =\exp \left( -2l\bar{s}\left( x_{0},r\right) \right) $ and
dividing by
\begin{equation}
\exp \left( 3l\bar{s}\left( x_{0},r\right) \right) +\eta ,
\end{equation}
we have
\begin{equation*}
\int\limits_{B\left( x_{0},qr\right) }G_{B\left( x_{0},2r\right)
}^{x_{0}}\alpha \left( u,u\right) \left( dx\right)
\end{equation*}
\begin{equation*}
+\frac{\left( e^{-2l\bar{s}\left( x_{0},r\right) }+cq^{\nu -2}e^{-l\gamma
\mu \left( x_{0},r\right) }\frac{c_{0}^{2}\left( x_{0}\right) }{c_{1}\left(
qr\right) }\right) }{e^{3l\bar{s}\left( x_{0},r\right) }+e^{-2l\bar{s}\left(
x_{0},r\right) }}\underset{B\left( x_{0},qr\right) }{\sup }\left| u-\bar{u}%
\right| ^{2}
\end{equation*}
\begin{equation*}
\leq \frac{\left[ c_{q}c_{1}^{4}\left( qr\right) e^{3l\bar{s}\left(
x_{0},r\right) }s_{0}^{18}\left( x_{0}\right) +1+2\eta cs_{0}^{8}\left(
x_{0}\right) c_{1}^{2}\left( r\right) e^{2l\bar{s}\left( x_{0},r\right) }%
\right] }{e^{5l\bar{s}\left( x_{0},r\right) }+1}\underset{B\left(
x_{0},r\right) }{\sup }\left| u-\bar{u}\right| ^{2}
\end{equation*}
\begin{equation*}
+\frac{1}{1+e^{-5l\bar{s}\left( x_{0},r\right) }}\int\limits_{B\left(
x_{0},r\right) }G_{B\left( x_{0},2r\right) }^{x_{0}}\alpha \left( u,u\right)
\left( dx\right)
\end{equation*}
\begin{equation*}
\leq \frac{3c_{q}c_{1}^{2}\left( qr\right) s_{0}^{18}\left( x_{0}\right)
e^{3l\bar{s}\left( x_{0},r\right) }}{e^{5l\bar{s}\left( x_{0},r\right) }+1}%
\underset{B\left( x_{0},r\right) }{\sup }\left| u-\bar{u}\right| ^{2}
\end{equation*}
\begin{equation}
+\frac{1}{1+e^{-5l\bar{s}\left( x_{0},r\right) }}\int\limits_{B\left(
x_{0},r\right) }G_{B\left( x_{0},2r\right) }^{x_{0}}\alpha \left( u,u\right)
\left( dx\right) .
\end{equation}
On the other hand the first member of previous inequality becomes
\begin{equation*}
\int\limits_{B\left( x_{0},qr\right) }G_{B\left( x_{0},2r\right)
}^{x_{0}}\alpha \left( u,u\right) \left( dx\right)
\end{equation*}
\begin{equation*}
+\frac{\left( 1+cq^{\nu -2}e^{2l\bar{s}\left( x_{0},r\right) }e^{-l\gamma
\mu \left( x_{0},r\right) }\frac{c_{0}^{2}\left( x_{0}\right) }{c_{1}\left(
qr\right) }\right) }{e^{5l\bar{s}\left( x_{0},r\right) }+1}\underset{B\left(
x_{0},qr\right) }{\sup }\left| u-\bar{u}\right| ^{2}.
\end{equation*}
This means that
\begin{equation*}
\int\limits_{B\left( x_{0},qr\right) }G_{B\left( x_{0},2r\right)
}^{x_{0}}\alpha \left( u,u\right) \left( dx\right) \leq \frac{%
3c_{q}c_{1}^{2}\left( qr\right) s_{0}^{18}\left( x_{0}\right) e^{3l\bar{s}%
\left( x_{0},r\right) }}{e^{5l\bar{s}\left( x_{0},r\right) }+1}\underset{%
B\left( x_{0},r\right) }{\sup }\left| u-\bar{u}\right| ^{2}
\end{equation*}
\begin{equation*}
+\frac{1}{1+e^{-5l\bar{s}\left( x_{0},r\right) }}\int\limits_{B\left(
x_{0},r\right) }G_{B\left( x_{0},2r\right) }^{x_{0}}\alpha \left( u,u\right)
\left( dx\right) \leq 2cc_{1}^{2}\left( qr\right) e^{-2l\bar{s}\left(
x_{0},r\right) }\underset{B\left( x_{0},r\right) }{\sup }\left| u-\bar{u}%
\right| ^{2}
\end{equation*}
\begin{equation}
+\frac{1}{1+e^{-5l\bar{s}\left( x_{0},r\right) }}\int\limits_{B\left(
x_{0},r\right) }G_{B\left( x_{0},2r\right) }^{x_{0}}\alpha \left( u,u\right)
\left( dx\right) .
\end{equation}
In the last inequality $\underset{B\left( x_{0},r\right) }{\sup }\left| u-%
\bar{u}\right| ^{2}$ can be related to the energy by \cite{Mosco}
\begin{equation*}
\underset{B\left( x_{0},r\right) }{\sup }\left| u-\bar{u}\right| ^{2}\leq
c\left( \underset{B\left( x_{0},R\right) }{osc}u\right) ^{2}\exp \left(
-\int\limits_{r}^{R}\exp \left( -2\gamma \mu \left( x_{0},\rho \right) \frac{%
d\rho }{\rho }\right) \right)
\end{equation*}
\begin{equation*}
\leq \frac{1}{m\left( B\left( x_{0},kR\right) \right) }\exp \left(
-\int\limits_{r}^{R}\exp \left( -2\gamma \mu \left( x_{0},\rho \right) \frac{%
d\rho }{\rho }\right) \right) \int\limits_{B\left( x_{0},kR\right) }\left| u-%
\bar{u}\right| ^{2}m\left( dx\right)
\end{equation*}
\begin{equation*}
\leq ce^{2l\gamma \mu \left( x,r\right) }c_{1}^{2}\left( R\right)
\int\limits_{B\left( x_{0},k^{2}R\right) }G_{B\left( x_{0},2q^{-1}r\right)
}^{x_{0}}\alpha \left( u,u\right) \left( dx\right) \cdot
\end{equation*}
\begin{equation}
\cdot \exp \left( -\int\limits_{r}^{R}\exp \left( -2\gamma \mu \left(
x_{0},\rho \right) \frac{d\rho }{\rho }\right) \right) .
\end{equation}
Then
\begin{equation*}
\int\limits_{B\left( x_{0},qr\right) }G_{B\left( x_{0},2r\right)
}^{x_{0}}\alpha \left( u,u\right) \left( dx\right)
\end{equation*}
\begin{equation*}
\leq cc_{1}^{2}\left( qr\right) c_{1}^{2}\left( R\right) \exp \left(
-\int\limits_{r}^{R}\exp \left( -2\gamma \mu \left( x_{0},\rho \right) \frac{%
d\rho }{\rho }\right) \right) \cdot
\end{equation*}
\begin{equation}
\cdot \int\limits_{B\left( x_{0},k^{2}R\right) }G_{B\left(
x_{0},2q^{-1}r\right) }^{x_{0}}\alpha \left( u,u\right) \left( dx\right) +%
\frac{1}{1+e^{-5l\bar{s}\left( x_{0},r\right) }}\int\limits_{B\left(
x_{0},r\right) }G_{B\left( x_{0},2q^{-1}r\right) }^{x_{0}}\alpha \left(
u,u\right) \left( dx\right) .
\end{equation}
If we consider $u\in \mathcal{D}\left[ a,B\left( x_{0},R_{0}\right) \right] $%
, $B\left( x_{0},R_{0}\right) \subset X$ and $R<kR$, one gets
\begin{equation}
\int\limits_{B\left( x_{0},r\right) }G_{B\left( x_{0},2q^{-1}r\right)
}^{x_{0}}\alpha \left( u,u\right) \left( dx\right) \leq \int\limits_{B\left(
x_{0},kR\right) }G_{B\left( x_{0},2q^{-1}r\right) }^{x_{0}}\alpha \left(
u,u\right) \left( dx\right)
\end{equation}
and
\begin{equation*}
\int\limits_{B\left( x_{0},qr\right) }G_{B\left( x_{0},2r\right)
}^{x_{0}}\alpha \left( u,u\right) \left( dx\right)
\end{equation*}
\begin{equation*}
\leq cc_{1}^{4}\left( qr\right) \exp \left( -\int\limits_{r}^{R}\exp \left(
-2\gamma \mu \left( x_{0},\rho \right) \frac{d\rho }{\rho }\right) \right)
\int\limits_{B\left( x_{0},k^{2}R\right) }G_{B\left( x_{0},2q^{-1}r\right)
}^{x_{0}}\alpha \left( u,u\right) \left( dx\right)
\end{equation*}
\begin{equation*}
+\frac{1}{e^{5l\bar{s}\left( x_{0},r\right) }+1}\int\limits_{B\left(
x_{0},kR\right) }G_{B\left( x_{0},2q^{-1}r\right) }^{x_{0}}\alpha \left(
u,u\right) \left( dx\right)
\end{equation*}
\begin{equation*}
\leq cc_{1}^{4}\left( qr\right) \exp \left( -\int\limits_{r}^{R}\exp \left(
-2\gamma \mu \left( x_{0},\rho \right) \frac{d\rho }{\rho }\right) \right)
\int\limits_{B\left( x_{0},k^{2}R\right) }G_{B\left( x_{0},2q^{-1}r\right)
}^{x_{0}}\alpha \left( u,u\right) \left( dx\right)
\end{equation*}
\begin{equation}
+\exp \left( -\int\limits_{r}^{R}\exp \left( -5l\gamma \mu \left( x_{0},\rho
\right) \right) \frac{d\rho }{\rho }\right) \int\limits_{B\left(
x_{0},k^{2}R\right) }G_{B\left( x_{0},2q^{-1}r\right) }^{x_{0}}\alpha \left(
u,u\right) \left( dx\right) .
\end{equation}
If we make the choice $\gamma \mu \left( x,r\right) \leq o\left( \log \log
\frac{1}{r}\right) $, we obtain by the previous inequality
\begin{equation*}
\psi \left( r\right) \leq 2\exp \left( -\int\limits_{r}^{R}\exp \left( -\log
\log \frac{1}{\rho }\right) \frac{d\rho }{\rho }\right) \psi \left(
k^{2}R\right)
\end{equation*}
\begin{equation}
\leq \left( \frac{\log \frac{1}{R}}{\log \frac{1}{r}}\right) \psi \left(
k^{2}R\right) ,
\end{equation}
where
\begin{equation}
\psi \left( r\right) =\int\limits_{B\left( x_{0},r\right) }G_{B\left(
x_{0},2q^{-1}r\right) }^{x_{0}}\alpha \left( u,u\right) \left( dx\right) .
\end{equation}

\end{document}